\documentclass{amsart}
\usepackage{amssymb,amsmath,amsthm, verbatim, latexsym}

\theoremstyle{plain}
\newtheorem{thrm}{Theorem}[section]
\newtheorem{lemma}[thrm]{Lemma}
\newtheorem{prop}[thrm]{Proposition}
\newtheorem{cor}[thrm]{Corollary}

\newtheorem{dfn}[thrm]{Definition}

\newtheorem{example}{\textbf{Example}}

\numberwithin{equation}{section}

\begin{document}

\newcommand{\Dxk}{\frac{\partial}{\partial x_k}}
\newcommand{\Dhk}{\frac{\partial}{\partial h_k}}
\newcommand{\Rm}{\mathbb R^m}
\newcommand{\Je}{J_\epsilon}
\newcommand{\Jei}{J_{\epsilon_i}}
\newcommand{\Kone}{K_\epsilon(x,h)}
\newcommand{\eh}{\epsilon h}
\newcommand{\IO}{\int_\Omega}
\newcommand{\Rn}{\mathbb R^n}
\newcommand{\Ia}{\int_{\alpha B}}
\newcommand{\Iaa}{\int_{(\alpha+1) B}}
\newcommand{\Om}{\Omega}
\newcommand{\BVX}{BV_{H}}
\newcommand{\Lone}{L^1(\Omega)}
\newcommand{\Hn}{\mathbb H^n}
\newcommand{\aB}{\alpha B}
\newcommand{\aBa}{(\alpha+1)B}
\newcommand{\dist}{\operatorname{dist}}
\newcommand{\loc}{\text{loc}}
\newcommand{\diam}{\operatorname{diam}}
\newcommand{\Var}{\operatorname{Var}}
\newcommand{\diver}{\operatorname{div}}
\newcommand{\iavg}{{\lim_{r \rightarrow 0} \int \!\!\!\!\!\! - }_{\!\!B(x,r)}}
\newcommand{\iavgr}{{\lim_{r \rightarrow 0} \int \!\!\!\!\!\! - }_{\!\!B_\rho(x,r)}}
\newcommand{\gr}{{\textrm{H}}}
\newcommand{\rb}{\partial _ {H}^*}
\newcommand{\mb}{\partial _ { *, H}}
\newcommand{\pme}{|\partial E | _ {H}}
\newcommand{\sdq}{\mathcal{S}^{Q-1} _d}
\newcommand{\sddq}{\mathcal{S}^{Q-1}_d}
\newcommand{\br}{B_\rho (x,r)}
\newcommand{\nh}{\nu_\Hs (x)}
\newcommand{\Hs}{\mathbb H^n}
\newcommand{\gro}{\mathbb{G}}
\title[Extension and Trace Theorem for $\gr$-BV]{An Extension and Trace Theorem for Functions
of $\gr$-Bounded
Variation in Carnot Groups of Step 2}
\author{Christina Selby}
\address{Mathematics Deparment\\
         Purdue University\\
	 West Lafayette, IN  47907}
\email{cselby@math.purdue.edu}
\date{\today}


\maketitle
\begin{abstract}
This paper provides an extension of a function $u \in \BVX(\Omega)$ to a function $u_0 \in
\BVX(\gro)$, when $\Omega \subset \gro$ is ``$\gr$-admissibile,'' and $\gro$ is a step 2 Carnot
Group. It is shown that $\gr$-admissible domains include
non-characteristic domains and domains in groups of Heisenberg type which have a
partial symmetry about characteristic points.  An example is given of a
domain $\Omega$ that is $C^{1, \alpha}$, $\alpha <1$, that is not $\gr$-admissible.  Further,
when $\Omega$ is $\gr$-admissible a trace
theorem is proved for $u \in \BVX(\Omega)$ . 
\end{abstract}

\section{Introduction}
There has been significant progress in the study of functions of $\gr$-bounded variation and
sets of finite $\gr$-perimeter in such works as \cite{AMB, AM, CDG1, FSSC2, FSSC1, GN, MAG}.  This paper utilizes these results to
further this progress.  In particular, the structure theorem for sets of finite $\gr$-perimeter
in step 2 Carnot groups, proved by Franchi, Serapioni, and Serra Cassano in \cite{FSSC1},
provides a way to define
an ``$\gr$-admissible domain.''  Their result states that for a set $E$ with finite $\gr$-perimeter
measure, the reduced boundary of $E$, $\rb E$, is $\gr$-rectifiable.  Further, the
$\gr$-perimeter measure of $E$ is the $(Q-1)$-spherical Hausdorff measure with respect to $d$
restricted to $\rb E$, where $d$ is a distance equivalent to the Carnot-Carath\'eodory
distance.  This measure is denoted $\sdq$. These $\gr$-admissible domains are significant, for if $\Omega$
is $\gr$-admissible and $u \in \BVX(\Omega)$, one may extend $u$ to be zero outside of $\Omega$
and obtain a function $u_0$ in $\BVX(\gro)$.  The $\gr$-admissibility condition given is an analogue to the condition stated in Meyer's and Ziemer's work in \cite{MZ}.  The
definition of $\gr$-admissible is given below, where $\mb E$ is the measure-theoretic
boundary of $E$:  
\begin{dfn} [$\gr$-Admissible Domain]
A bounded domain $\Omega$ of finite $\gr$-perimeter is said to be $\gr$-admissible if
the following two conditions are satisf{}ied:
\begin{itemize}
\item[(i)] $\sdq (\partial \Omega \backslash \mb \Omega) =0$.
\item[(ii)] There is a constant $M=M(\Omega)$ such that for each $x \in \partial
\Omega$ there is a $d$-ball $B(x,r)$ with
\[
\sdq (\mb E \cap \mb \Omega ) \leq M \ \sdq (\mb E \cap \Omega )
\]
for all $E \subset \overline{\Omega} \cap B(x,r)$.
\end{itemize}
\end{dfn}
\noindent Using Meyer's and Ziemer's Euclidean definition, Lipschitz domains can be shown to be admissible quite easily
using the Gauss-Green theorem.  However, a wide class of $\gr$-admissible domains in the setting of step 2 Carnot Groups is not so easily obtained.  This paper
shows that $\gr$-admissible domains include
non-characteristic domains and domains in groups of Heisenberg type which have a
partial symmetry about characteristic points.  In proving condition (i) of $\gr$-admissibility a fact proved by Magnani in \cite{MAG} is
vital.  It
states the measure of the
characterstic set of a $C^1$ domain in a step 2 Carnot group has $\sdq$ measure zero.   This result is an extension of a result by Balogh for the Heisenberg group \cite{BA}. This paper
also gives an
example of a $C^{1,\alpha}$ domain, $\alpha < 1$ that is not
$\gr$-admissible.     

The extended function $u_0$ is not defined on $\partial \Omega$, but this paper defines the trace
of $u$ on $\partial \Omega$, $u^*$, as the sum of the upper and lower approximate limit
of $u_0$, $\mu$ and $\lambda$, respectively.  Here, the definition of upper and lower approximate limit is that of Federer \cite{FE}.  The set $F$, where the lower approximate
limit is strictly less than the upper
approximate limit is shown to be $\gr$-rectifiable.  Ambrosio and Magnani use a stronger definition for approximate continuity in \cite{AM}, and prove that the approximate discontinuity set is
$\gr$-rectifiable. This paper shows that the set $F$, however, has some additional useful properties. The stated rectifiability result, along with an implicit function theorem proved in \cite{FSSC1} lead to the following trace theorem:
\begin{thrm}
If $\Omega$ is an $\gr$-admissible domain, there is a constant $M = M(\Omega)$
such that
\[
\int _{\gro} |u^*| \ \ d |\partial \Omega|_\gr = C(n+m) \ \int _{\mb \Omega} |u^*| \
\ d\sdq \leq  M \ \parallel u \parallel
_{\BVX (\Omega)}
\]
whenever $u \in \BVX (\Omega).$
\end{thrm}

The additional useful properties of the set $F$ are stated below: 
\begin{thrm} 
Let $u \in \BVX (\gro)$, then
\begin{itemize}
\item[(i)] $F$ is countably $\gr$-rectif{}iable,
\item[(ii)] $- \infty < \lambda (x) \leq \mu (x) < \infty$ for $\sdq$ a.e. $x \in \gro$,
\item[(iii)] for $\sdq$ a.e. $y \in F$, there is a  vector
$\nu(y)$ such that {} $\nu_{A_s}(y) = \nu(y)$ whenever $\lambda (y) < s < \mu
(y)$.
\item[(iv)] For all $y$ as in (iii), with $-\infty < \lambda (y) < \mu (y) <
\infty$ there are Lebesgue measurable sets $F^+$ and $F^-$ such that
\begin{align*}
&\lim_{r \rightarrow 0} \frac{|F^- \cap \tau _y (S_\gr ^- ( \nu  (y))) \cap
B(y,r)|}{|\tau_y(S_\gr ^- ( \nu (y))) \cap B(y,r)|}\\
= &\lim_{r \rightarrow 0} \frac{|F^+ \cap \tau _y (S_\gr ^+( \nu  (y))) \cap
B(y,r)|}{|\tau _y (S_\gr ^+ ( \nu (y))) \cap B(y,r)|} =1
\end{align*}
and
\[
\lim_{\substack{x \rightarrow y \\ x \in F^- \cap \tau _y (S_\gr ^- ( \nu (y)))}}
u(x) = \mu (y), \qquad \lim_{\substack{x \rightarrow y \\ x \in F^+ \cap \tau _y
(S_\gr ^+ (\nu(y)))}}
u(x) = \lambda (y).
\]
\end{itemize}
\end{thrm}
In (iv), the sets $S_\gr ^+ ( \nu  (y))$ and $S_\gr ^- ( \nu  (y))$ are the
half-spaces obtained in the blow-up theorem found in \cite{FSSC1}.  A Poincar\'e inequality proved by Garofalo and Nheiu in \cite{GN}
for functions in $\BVX(\br)$, $\br$ a gauge ball, along with part (iii) are used to prove some integral properties of a function $u \in
\BVX(\gro)$. These are stated below, where $U(x)$ is the average of the upper and lower
approximate limit of $u$:
\begin{thrm}
Assume $u \in \BVX (\gro)$.  Then 
\begin{itemize}
\item[(i)] ${\lim_{r \rightarrow 0} \int \!\!\!\!\! - }_{\br}|u-U(x)| ^{\frac{Q}{Q-1}} \,dh =0$  for 
$\sdq$  a.e.  $x \in
\gro\backslash F$, and
\item[(ii)]for $\sdq$ a.e. $x \in F$, there exists a  vector $\nu =
\nu(x)$ such that
\[
\lim_{r \to 0}  {\int \!\!\!\!\!\! - }_{\!\!\br \cap S_\gr^- (\nu)} |u - \mu (x)|
^\frac{Q}{Q-1} \, dh = 0,
\]
and
\[
\lim_{r \to 0}  {\int \!\!\!\!\!\! - }_{\!\!\br \cap S_\gr^+ (\nu)} |u - \lambda (x)|
^\frac{Q}{Q-1} \, dh =0.
\]
\end{itemize}
\end{thrm}
  
Using the definition of Ambrosio, Fusco, and Pallara in \cite{A}, for approximate jump point, one can then observe
that a function $u \in \BVX(\gro)$ has an approximate jump point for $\sdq$ a.e. point in $F$.   
  
\section{Notation and Preliminary Results}
A Carnot group $\gro$ of step 2 is a simply connected Lie group whose Lie
algebra $\mathfrak{g}$ admits a step 2 stratification.  This means there are
subspaces $V_1, V_2$ such that
\[
\mathfrak{g} = V_1 \oplus V_2 , \qquad [V_1, V_1] = V_2, \qquad [V_1, V_2]=0,
\]
where $[V_i, V_k]$ is the subspace of $\mathfrak{g}$ generated by the commutators
$[X,Y]$, where $X \in V_i$, $Y \in V_k$.  Throughout this paper, let $m =
\textrm{dim}(V_1)$, and $n=\textrm{dim}(V_2)$.  Let $\{e_1, \ldots, e_m\}$ denote an
orthornormal basis for $V_1$, and $\{\epsilon_1, \ldots, \epsilon_n\}$ denote an
orthonormal basis for $V_2$.  Let $X= \{X_1, \ldots, X_m\}$ be the family of
left-invariant vector fields where $X_i(0)=e_i$.  The set $\{X_1, \dots, X_m\}$
along with all of its commutators generates $\mathfrak{g}$.  The exponential map
is known to be a global diffeomorphism, and using it, one is able to identify an
element of $\gro$ with an element of $\mathbb{R}^{m+n}$.  In particular, for $p
\in \gro$, $p=\textrm{exp}(p_1X_1+\cdots+p_{m+m}X_{m+n})$, and can be identified
with $(p_1, \ldots, p_{m+n}) \in \mathbb{R}^{m+n}$.  Thus $\gro$ can be
identified with $(\mathbb{R}^{m+n},\cdot)$, where the group operation $\cdot$ is
determined by the Campbell-Hausdorff formula. 

The \emph{horizontal bundle} $H\gro$ is the subbundle of the tangent bundle
$T\gro$ that is spanned by $X_1,\ldots, X_m$.  The fibers of $H\gro$ are
\[
H\gro_x= \textrm{span} \{X_1(x), \ldots, X_m(x)\}, \qquad x \in \gro.
\] 
Each fiber of $H\gro$ is endowed with a scalar product $\langle \cdot, \cdot
\rangle_x$ and a norm $|\cdot|_x$ that makes the basis $\{ X_1(x), \ldots,
X_m(x)\}$ an orthonormal basis.  So for $p=\sum_{i+1}^mp_iX_i(x)=(p_1,\ldots,
p_m)$, and $q=\sum_{i+1}^mq_iX_i(x)=(q_1,\ldots,
q_m)$, $\langle p, q \rangle = \sum_{i=1}^m p_iq_i$, and $|p|_x^2=\langle p,p
\rangle_x$.  
The sections of $H\gro$ are called \emph{horizontal sections} and a vector in
$H\gro$ is called a \emph{horizontal vector}.  A horizontal section $\phi$
can be identified with a function $\phi = (\phi_1, \ldots,
\phi_m\}:\mathbb{R}^{m+n} \rightarrow \mathbb{R}^m$ by identifying it with its
canonical coordinates with respect to $\{X_1(x), \ldots, X_m(x)\}$.  The
notation $\langle \phi, \psi \rangle$ represents $\langle \phi(x), \psi(x)
\rangle_x$.  The \emph{horizontal divergence} of $\phi: \mathbb{R}^{m+n} \rightarrow
\mathbb{R}^m$ is defined as
\[
\textrm{div}_\gr \phi (x) := \sum_{i=1}^m X_i\phi_i(x),
\]
and the \emph{horizontal gradient} of $u \in C^1(\gro)$ is defined as 
\[
\nabla_{\gr}u(x) := (X_1u (x), \ldots, X_mu(x)).
\]
Translations and dilations are defined on $\gro$ as follows:
\[
\tau_p(q) = p\cdot q, \qquad \delta_r(p) = (rp_1, \dots, rp_m, r^2p_{m+1},
\ldots, r^2p_{m+n}).
\]
Then $Q:= m + 2n$ is called the \emph{homogeneous dimension} of $\gro$. 

The results of Franchi, Serapioni, and Serra Cassano utilize the distance
defined below:
\[
d(x,y) = d(y^{-1} \cdot x, 0),
\]
where for $p=(p_1, \ldots, p_{m+n}) \in \mathbb{R}^{m+n}$,
\begin{equation}\label{dis}
d(p,0)= \textrm{max}\{ ||(p_1, \ldots, p_m)||_{\mathbb{R}^m},
\epsilon ||(p_{m+1},\ldots, p_{m+n})||_{\mathbb{R}^n}^{1/2}\}.
\end{equation}
Here $\epsilon \in (0,1)$ is a constant depending on the group structure. 
Throughout this paper, $B(x,r):=\{ y \in \gro : d(x,y) < r\}$.  This distance $d$
is
known to be equivalent to the Carnot Carath\'eodory distance, and thus also
equivalent to the gauge distance (see \cite{FSSC1}). The gauge distance is given
below:
\[
\rho(x,y)=\rho(y^{-1}\cdot x, 0),
\]
\[
\rho(x, 0) = ( (x_1^2 + \cdots + x_m^2)^2 + x_{m+1}^2 + \cdots + x_{m+n}^2 )
^{1/4}.
\]
The gauge balls $B_\rho(x,r) := \{y \in \gro : \rho(x,y) < r\}$ are known to be
\emph{PS-domains}, see \cite{GN}, for which the isoperimetric inequality and a Poincar\'e
inequality hold.  These
inequalities will be stated shortly.  The $(Q-1)$-dimensional spherical Hausdorff
measure with respect to $d$ will be denoted $\sdq$.  Proposition 2.4 in
\cite{FSSC1} states that the diameter of $B(x,r)$ with respect to $d$ is $2r$.  

The $\gr$-variation of $u\in L^1_{loc}(\Omega)$ with respect to $\Omega$ is given as
\[
Var_\gr (u; \Omega) = \sup \{ \int_\Omega u \,div_H \phi \,dh : \phi \in C_0 ^1 (\Omega ;
\mathbb{R}^{m}),
|\phi(P)|\leq 1 \}, 
\]
where $dh$ is $(m+n)$-dimensional Lebesgue measure on $\mathbb{R}^{m+n}$  The
function $u$ is said to belong to $\BVX(\Omega)$ if $Var_\gr(u; \Omega) < \infty$. 
The $\gr$-perimeter measure of
$\Omega$ with respect to the set $E$ is given as $\pme(\Omega) := \Var_\gr (\chi_E; \Omega)$.
The set $E$ is
said to have finite $\gr$-perimeter, or be $\gr$-Caccioppoli, if $\pme(\gro) < \infty$.  By Riesz representation theorem $\pme$ is a Radon measure on $\gro$ and there
exists a measurable section $\nu_E$ of $H\gro$ such that
\[
\int_E \textrm{div}_\gr \phi \, dh = -\int_\gro \langle \phi, \nu_E \rangle \, d\pme,
\]
for any $\phi \in C_0^1(\gro, H\gro)$.  The section $\nu_E$ is called the \emph{generalized
inward normal} to $E$.  
The following proposition will be
used in a following section and is proved in \cite{CDG1}.
\begin{prop}\label{form}
If $E$ is $\gr$-Caccioppoli with Euclidean $C^1$ boundary, then 
\[
\pme(\Omega) = \int_{\partial E \cap \Omega} \big ( \sum_{i=1}^m \langle X_i, n
\rangle_{\mathbb{R}^m}^2\big)^{1/2} \, d \mathcal{H}^{m+n-1},
\]
where $\mathcal{H}^{m+n-1}$ is Euclidean $(m+n-1)$-dimensional Hausdorff measure and $n=n(x)$
is the Euclidean unit outward normal to $\partial E$.  
\end{prop}

If $\Omega$ can be described as $\Omega= \{\phi < 0 \}$, then $\phi$ is called a \emph{defining
function} for $\Omega$.  Further, $n= -\nabla \phi / |\nabla \phi|$.  The above proposition
will be used with such a defining function.  

The following three theorems are proved in \cite{GN}.  Theorem \ref{coarea} was also proved independently in \cite{FSSC3}.
\begin{thrm}\label{iso}
There is a positive constant $c>0$ such that for any $\gr$-Caccioppoli set $E$, for all $x\in
\gro$, and $r>0$,
\[
\textrm{min}\{|E\cap B_\rho(x,r)|, |E^c \cap B_\rho(x,r)|\} ^{\frac{Q-1}{Q}} \leq
c\pme(B_\rho(x,r)), \textrm{and}
\]
\[
\textrm{min}\{|E|, |E^c|\}^\frac{Q}{Q-1} \leq c \pme(\gro),
\]
where $|E|$ is the $(m+n)$-dimensional Lebesgue measure of $E$.
\end{thrm}

\begin{thrm}\label{poin}
For any $u \in \BVX (\br)$, one has $C=C(m+n)$ such that
\[
\|u - u_{x,r}\|_{L^\frac{Q}{Q-1}(\br)} \leq C \,Var_\gr(u; \Omega),
\]
where $u_{x,r}$ denotes the integral average of $u$ over $\br$.
\end{thrm}

\begin{thrm}{(Coarea Formula)}\label{coarea}
Let $u \in \BVX(\Omega)$.  Then
\[
Var_{\gr} (u ; \Omega) = \int_{-\infty}^\infty |\partial A_t|_\gr(\Omega) \, dt,
\]
where $A_t= \{ x \in \Omega : u(x) > t\}$. 
\end{thrm}
Now for some measure-theoretic defintions.
\begin{dfn}[Reduced Boundary]
Let $E$ be an $\gr$-Caccioppoli set; then $x \in \rb E$, the reduced boundary of $E$, if 
\begin{gather*}
\pme (B(x,r)) > 0 \qquad \textrm{for any } r > 0, \\
\textrm{there exist} \quad \iavg \nu_E \: d\pme, \quad \textrm{and} \\
\iavg \nu_E \: d\pme = 1.
\end{gather*}
\end{dfn}

\begin{dfn}[Measure-Theoretic Boundary]\label{MTB}
Let $E$ $\subset \gro$ be a measurable set, then $x \in \mb E$, the measure theoretic boundary of $E$, if
\[
\limsup_{ r \rightarrow 0^+} \frac{|E \cap B(x,r)|}{|B(x,r)|} > 0 \qquad and \qquad \limsup_{ r
\rightarrow 0^+} \frac {|E^c \cap B(x,r)|}{|B(x,r)|} >0.
\]
\end{dfn}

The following lemma and theorem can be found in \cite{FSSC1}:
\begin{lemma} \label{mea}
Let $E$ be an $\gr$-Caccioppoli set, then
\begin{itemize}
\item[(i)]$\rb E \subseteq \mb E  \subseteq \partial E$
\item[(ii)] $\sdq (\mb E \ \backslash \ \rb E) =0.$
\end{itemize}
\end{lemma}

\begin{thrm}[Gauss-Green Theorem]\label{ggt}
Let $E$ be an $\gr$-Caccioppoli set, then
\[
-\int_E \diver_\gr \phi \, dh = \theta_d \int_{\mb E} \langle \nu_E, \phi \rangle  \, d\sdq, \qquad
\textrm{for all } \phi \in C_0^1(\gro, H\gro).
\]
\end{thrm}

The following definition, found in \cite{FSSC1} is needed for the statement of the implicit
function theorem. Here $C_\gr^1(\Omega)$ is the collection of functions $u$ with 
distributional derivatives $X_iu$ that are continuous in $\Omega$. 
\begin{dfn}
$S\subset \gro$ is an $\gr$-regular hypersurface if for every $x \in S$ there exists a
neighborhood $\mathcal{U}$ of $x$ and a function $f \in C_\gr^1(\mathcal{U})$ such that
\[
S \cap \mathcal{U} = \{ y \in \mathcal{U}: f(y)=0\}, \textrm{and}
\]
\[
\nabla_\gr f(y) \neq 0 \qquad \textrm{for } y \in \mathcal{U}.
\]
\end{dfn}
\begin{thrm}[Implicit Function Theorem]\label{IFT}
Let $\Omega$ be an open set in $\mathbb{R}^{m+n}$, $0 \in \Omega$ and let $f \in C_ \gr ^1
(\Omega)$ be such that $X_1 f(0) >0, f(0) = 0$.  Define
\[
E=\{x \in \Omega : f(x) <0\}, \quad S= \{x\in \Omega : f(x) =0\},
\]
and, for $\delta >0$, $h>0$
\[
I_\delta = \{\xi=(\xi_2, \ldots, \xi_{m+n}) \in \mathbb{R}^{m+n-1}, |\xi_j|\leq \delta\},
\qquad J_h = [h, h].
\]
If $\xi = (\xi_2, \ldots, \xi_{m+n}) \in \mathbb{R}^{m+n-1}$ and $t \in J_h$, denote by
$\gamma(t, \xi)$ the integral curve of the vector field $X_1$ at the time $t$ issued from
$(0,\xi) \in \mathbb{R}^{m+n}$, i.~e.
\[
\gamma(t, \xi)= \textrm{exp}(tX_1)(0, \xi).
\]
Then there exists $\delta$, $h >0$ such that the map $(t, \xi) \rightarrow \gamma(t, \xi)$
is a homeomorphism of a neighborhood of $J_h \times I_\delta$ onto an open subset of
$\mathbb{R}^{m+n}$, and, if $\mathcal{U} \subset \subset \Omega$ is the image of
$\textrm{Int}(J_h\times I_\delta)$ through this map, 
\begin{itemize}
\item[(i)]$E$ has f{}inite $\gr$ -perimeter in $\mathcal{U}$;
\item[(ii)] $\partial E \cap \Omega = S \cap \mathcal{U}$;
\item[(iii)] $\nu _E (s) = - \nabla_\gr f(x) / |\nabla _ \gr f(x)|_ x$ for all $x \in S
\cap \mathcal{U}$,
where $\nu _E$ is the generalized inner unit normal.  Moreover, there exists a unique
function 
\[
\phi = \phi(\xi): I_\delta \rightarrow J_h
\]
such that the following parameterization holds:  if $\xi \in I_\delta$ and we put $\Phi(\xi)
= \gamma(\phi(\xi), \xi)$, then 
\[
S \cap \tilde{\Omega} = \{ x \in \tilde{\mathcal{U}} : x = \Phi(\xi), \xi \in I_\delta \};
\]
\[
\phi \textrm{ is continuous};
\]
the $\gr$-perimeter has an integral representation
\[
\pme(\tilde{\mathcal{U}}) = \int_{I_\delta} \frac{\sqrt{\sum_{i=1}^m |X_j
f(\Phi(\xi))|^2}}{X_1f(\Phi(\xi))} \, d\xi.
\]
\end{itemize}
\end{thrm}
\medskip
\begin{thrm}[Structure Theorem]\label{ST}
If $E \subseteq \gro$ is an $\gr$-Caccioppoli set, then
\[
\rb E \textrm{ is } (Q-1) \textrm{-dimensional }\gr \textrm{-rectif{}iable,}
\]
that is, $\rb E = N \cup \bigcup _{h=1} ^{\infty} K_h$, where $\mathcal{H} _d
^{Q-1} (N) =0$
and $K_h$ is a compact subset of a $\gr$-regular hypersurface $S_h$;
\[
\nu _E (p) \textrm{ is } \gr \textrm{-normal to } S_h \textrm{ at } p, \qquad \forall p \in K_h,
\]
\[
\pme =  \theta_d \sdq \lfloor \rb E,
\]
where
\[
\theta_d
=\frac{\omega_{m-1}\omega_{n}\epsilon^n}{\omega_{Q-1}}=\frac{1}{\omega_{Q-1}}\mathcal{H}^{m+n-1}(\partial
S_\gr^+(\nu_E(0)\cap B(0,1)).
\]
Here $\epsilon$ is as in \ref{dis}.
\end{thrm}
Note that by lemma \ref{mea}, the measure theoretic boundary of $E$ is also $\gr$-rectifiable.

The half-spaces $S_\gr ^+(\nu _E (P))$ and $S_\gr ^-(\nu _E (P))$ appear in the blow-up
theorem stated in theorem 3.1 of \cite{FSSC1}.  The hyper-planes can be thought of as
``approximate tangent planes'' to the boundary of $E$ at $p$. They are defined below:
\medskip
\[
S_ \gr ^+(\nu _E (p)) := \{q: \langle \pi _{p} q, \nu _E (p) \rangle_{p} \geq 0
\},
\]
\[
S_ \gr^- (\nu _E (p)) := \{ q: \langle \pi _ {p} q, \nu _E (p) \rangle_{p} \leq
0 \},
\]
where for $q=(q_1, \ldots, q_m, q_{m+1}, \dots, q_{m+n})$, 
\[
\pi _{p}q= \sum_{i=1}^m
q_i\,X_i(p)=(q_1, \dots, q_m, 0, \ldots, 0) \in \mathbb{R}^{m+n}.
\]

The reader may have noticed that the set where the ``measure-theoretic normal'' exists has
not been explicitly defined.  This set is often defined in classic geometric measure theory
texts.  There is no
need to distinguish this set in this work because of lemma \ref{mea} and the
following lemma which is lemma 3.3 in \cite{FSSC1}:
\begin{lemma}\label{fss}
Let $p \in \rb E$.  Then
\begin{gather*}
\lim_{r \rightarrow 0} \frac{| B(p,r) \cap E \cap \tau _p (S_\gr ^-(\nu _E
(p)))|}{|B(p,r)|} =0,\\
\lim_{r \rightarrow 0} \frac{| B(p,r) \cap E^c \cap \tau _p (S_\gr ^-(\nu _E
(p)))|}{|B(p,r)|} =0, \quad \textrm{and}\\
\lim_{r \rightarrow 0} \frac{ \pme (B(p,r))}{r^{Q-1}} = |\partial
S_\gr^+(\nu_E(p))|_\gr(B(0,1)).
\end{gather*}
\end{lemma}
The right-hand side of the last statement of the lemma is actually a constant, as can be seen
in theorem 3.1 of \cite{FSSC1}.

\section{Examples of $\gr$-admissible domains}
In the following sections, a function $u \in \BVX(\Omega)$ is shown to have an extension in $\BVX(\gro)$ when $\Omega$ is $\gr$-admissible.  This section provides examples of
$\gr$-admissible domains, and proves that if $\Omega$ is a $C^{1, \alpha}$ domain, $\alpha <1$, then it is not necessarily $\gr$-admissible.  The definition of an $\gr$-admissible domain
follows:

\begin{dfn} [$\gr$-Admissible Domain]
A bounded domain $\Omega$ of finite $\gr$-perimeter is said to be $\gr$-admissible if
the following two conditions are satisf{}ied:
\begin{itemize}
\item[(i)] $\sdq (\partial \Omega \backslash \mb \Omega) =0$.
\item[(ii)] There is a constant $M=M(\Omega)$ such that for each $x \in \partial
\Omega$ there is a ball $B(x,r)$ with
\[
\sdq (\mb E \cap \mb \Omega ) \leq M \ \sdq (\mb E \cap \Omega )
\]
for all $E \subset \overline{\Omega} \cap B(x,r)$.
\end{itemize}
\end{dfn}

Before stating any theorems about $\gr$-admissible domains, a few observations should be
made. In the classical setting, it is proved that Lipschitz domains are admissible using the
Gauss-Green theorem by setting $V=(1,0)$ and utilizing the nice form of the normal vector of a
graph (see Remark  5.10.2 in \cite{Z}).  Observe that if $\Omega$ is a $C^1$ domain with no characteristic points then without loss of generality,
for $x \in \partial \Omega$, there is a ball $B(x,r)$ such that $\nu_{\Omega,1} > K > 0$ in
$\partial \Omega \cap B(x,r)$.  Taking $V=(1,0)$, one obtains
\[
\int _{\mb E \cap \partial_{*, \gr} \Omega} \langle \nu _\Omega, V \rangle \, d\sddq =   \int
_{\mb E \cap \Omega} \langle \nu _E , V \rangle \, d\sddq,
\]
for $E \subset \overline{\Omega} \cap B(x,r)$, which implies
\[
K \sddq(\mb E \cap \mb \Omega ) \leq  \sddq (\mb E \cap \Omega ).
\]
 By the compactness of $\partial
\Omega$ a covering argument gives a uniform $K$ such that the above holds. Since $\Omega$ is $C^1$, it has finite $\gr$-perimeter. The following proposition, found in \cite{MAG}
provides condition (i) of $\gr$-admissibility for $C^1$ domains:
\begin{prop}\label{mag}
For $\Omega \subset \gro$ a $C^1$ domain, the characteristic set of $\Omega$ is
$\sdq$-negligible. 
\end{prop}
Since a $C^1$ domain has a horizontal unit normal at all points except characterstic
points, (i) holds.  This proves the following proposition.
\medskip
\begin{prop}
A bounded $C^1$ domain $\Omega \subset \gro$, with no characteristic points is $\gr$-admissible.
\end{prop}
\medskip
Since $\gr$-admissibility is a local condition, the proof of the previous proposition reveals
that  condition (ii) of $\gr$-admissibility will hold true for $C^1$
domains at points outside a neighborhood
of each characteristic point. Further, proposition \ref{mag} provides condition (i). Thus, to
prove a domain is $\gr$-admissible it is sufficient to show there is $M(\Omega)>0$ such that for
every characteristic point $x \in \partial \Omega$, there is $r>0$ such that 
\[
\sdq (\mb E \cap \mb \Omega ) \leq M \ \sdq (\mb E \cap \Omega )
\]
for all $E \subset \overline{\Omega} \cap B(x,r)$. The next proposition gives a condition for domains in a group of Heisenberg type to
be $\gr$-admissible when the domain a partial symmetry about characteristic points.  First, some definitions are needed.  Recall from section 2 that $\{ e_1, \ldots,
e_m\}$ and $\{\epsilon_1, \ldots, \epsilon_n\}$ denote orthonormal bases for $V_1$ and $V_2$ respectively.  
\begin{dfn}\label{kapmap}
In a Carnot group $\gro$ of step 2 with Lie algebra $\mathfrak{g} = V_1 \oplus V_2$,  the linear mapping $J: V_2 \to End(V_1)$ is defined by
\[
\langle J(\eta)\xi', \xi''\rangle = \langle [ \xi', \xi'' ], \eta \rangle, \qquad \eta \in V_2, \, \, \xi', \xi'' \in V_1. 
\]
\end{dfn}

\begin{dfn}\label{heis}
A Carnot group $\gro$ of step 2 is called of Heisenberg type if for every $\eta \in V_2$, such that $|\eta|=1$, the map $J(\eta): V_1 \to V_1$ is orthogonal.
\end{dfn}
From the definitions, the following observations can be made for a Carnot group of Heisenberg type:
\begin{equation}\label{mods}
|J(\eta)\xi| = |\eta||\xi|, \qquad \eta \in V_2, \xi \in V_1, 
\end{equation}
\begin{equation}\label{zero}
\langle J(\eta)\xi, \xi \rangle =0, \qquad \eta \in V_2, \xi \in V_1,
\end{equation}
\begin{equation}\label{ht}
\langle J(\eta)\xi , J(\eta')\xi \rangle = \langle \eta, \eta' \rangle |\xi|^2, \qquad \eta, \eta' \in V_2, \xi \in V_1.
\end{equation}

Formulas for the vector fields $X_i$, $i=1, \ldots , m$ have been obtained in \cite{DGN1} and \cite{FSSC1}.  They are as follows:
\[
X_i = \frac{\partial}{\partial x_i} + \frac{1}{2} \sum_{j=1}^m\sum_{l=1}^n b_{ji}^l x_j \frac{\partial}{\partial y_l}, \qquad i= 1,\ldots, m,
\]
where $b_{ij}^l = \langle [e_i, e_j] , \epsilon_l \rangle$.   
Using these formulas, the horizontal gradient of a
function $\phi \in C^1$ can be written using the mapping $J$.  Observe
\begin{align*}
X_i \phi &= \frac{\partial \phi}{\partial x_i} + \frac{1}{2} \sum_{j=1}^m\sum_{l=1}^n b_{ji}^l
x_j\frac{\partial \phi}{\partial y_l}\\
&=\frac{\partial \phi}{\partial x_i} + \frac{1}{2}\sum_{j=1}^m\sum_{l=1}^n \langle [e_j, e_i], \epsilon_l
\rangle x_j \frac{\partial \phi}{\partial y_l}\\
&=\frac{\partial \phi}{\partial x_i}+ \frac{1}{2} \langle [ \xi, e_i], \eta \rangle \\
&=\frac{\partial \phi}{\partial x_i} + \frac{1}{2} \langle J(\eta)\xi, e_i \rangle,
\end{align*}
where $\xi = \sum_{i=1}^m x_i e_i \in V_1$ and $\eta =
\sum_{l=1}^{n}\frac{\partial \phi}{\partial y_l} \epsilon_l \in V_2$. Therefore,
\begin{equation}
\nabla_H \phi = D_\xi \phi + \frac{1}{2} J(\eta) \xi,
\end{equation}
where $D_\xi$ is the standard Euclidean gradient in $\mathbb{R}^m$.  Thus, using \ref{mods}
\begin{equation}\label{nmod}
|\nabla_H \phi|^2 = |D_\xi|^2 + \langle D_\xi, J(\eta)\xi\rangle + \frac{1}{4}|\eta|^2|\xi|^2.
\end{equation}
So if $\phi$ is the defining function for a domain $\Omega$, then
\begin{equation}\label{nvector}
\nu_\Omega = -\frac{D_\xi \phi + \frac{1}{2}J(\eta)\xi}{\sqrt{|D_\xi \phi|^2 + \langle D_\xi \phi, J(\eta)\xi\rangle + \frac{1}{4}|\eta|^2|\xi|^2}}.
\end{equation}

A notion of a domain being cylindrically symmetric about characteristic points is defined in \cite{CG}.  A similar notion is defined in this paper:

\begin{dfn}
Let $\Omega \subset \gro$ be a bounded, connected, $C^1$ domain. Assume for each characteristic point $P$,
after a group translation that sends $P$ to the identity $e =0$, one can find a neighborhood $U$ of $e$ such that
\begin{itemize}
\item[(i)]$\partial \tilde{\Omega} \cap U = \{ y_n= -g(s, y_1, \ldots, y_{n-1}) \} \cap U$,
\item[(ii)] there exists $M(\Omega)>0$ such that $\left| \frac{\partial g}{\partial s}\right| \leq M|s|$ in $U$,
\item[(iii)]$g(0, y_1, \ldots, y_{n-1}) = \frac{\partial g}{\partial s}(0, y_1, \ldots, y_{n-1})
=0$, and
\item[(iv)]$g$ is $C^1$ with respect to the variables $y_1, \ldots, y_{n-1}$.
\end{itemize}
Here, $\tilde{\Omega} = P^{-1} \Omega$, $s=\sqrt{x_1^2 + \cdots + x_m^2} = |\xi|$. Then $\Omega$ is said to have partial symmetry near its characteristic
set.
\end{dfn}

\begin{thrm}
Assume that $\Omega \subset \gro$ has partial symmetry near its characterstic set.  Then $\Omega$ is $\gr$-admissible.
\end{thrm}
\begin{proof}[Proof:]
Note that the $\gr$-admissibility condition is invariant under group translation. So without loss of generality, consider a domain $\Omega$ with characteristic point at $e=0 \in
\mathbb{R}^m$.
Taking $\phi = y_n + g(s, y_1, \ldots, y_{n-1})$ as the defining function in the neighborhood $U$ of $0$ and using \ref{nvector} one computes
\[
\nu _\Omega  = -\frac{ \frac{\xi}{s}\frac{\partial g}{\partial s} +\frac{1}{2}J(\eta)\xi}{\sqrt{\left( \frac{\partial g}{\partial s} \right) ^2 + \frac{1}{4}|\eta|^2{s^2}}}.
\]
Here $\eta = \sum_{l=1}^{n-1}\frac{\partial g}{\partial y_l}\epsilon_l + \epsilon_n \in V_2$.   
\newline 
Let $B(0,r) \subset U$, and $E \subset B(0,r) \cap \overline{\Omega}$.  
Consider the vector
field $V=-\frac{J(\eta)\xi}{s}$.  This vector field is not $C_0^1(\gro ; H\gro)$, and so it can not be
used directly in the Gauss-Green theorem.  However, the theorem can still be used to show
$\gr$-admissibility.  Define
\[
V_\delta  = \begin{cases}
-\frac{J(\epsilon_n)\xi}{s}& \text{if $s>\delta$},\\
-\frac{J(\epsilon_n)\xi}{\delta}& \text{if $s \leq \delta$}.
\end{cases}
\]
$V_\delta$ is continuous, but not $C_0^1(\gro, H\gro)$.  Let $K: \mathbb{R}^{m} \to \mathbb{R}$,
$K(\xi) = K^*(s)$, be a standard mollifier and define
\[
V_{\delta,\tau}  = \begin{cases}
-\sum_{i=1}^m \frac{J(\epsilon_n)x_i * K_\tau}{s}e_i& \text{if $s>\delta$},\\
-\sum_{i=1}^m \frac{J(\epsilon_n)x_i * K_\tau}{\delta}e_i& \text{if $s \leq \delta$},
\end{cases}
\]
where $*$ is the standard convolution operator.
Then $V_{\delta, \tau}$ converges uniformly to $V_{\delta}$ on compact subsets of $\gro$ as 
$\tau$ goes to $0$.  Further, multiplying $V_{\delta, \tau}$ by a cutoff function which is
identically $1$ in a neighborhood of $E$ gives a vector field in $C_0^1(\gro; H\gro)$.  Call this
vector field $V_{\delta, \tau}$ as well. \newline
\noindent \textbf{Claim:} $\textrm{div}_{H} V_{\delta, \tau} =0$ on $E$.  First observe that since $V_{\delta,
\tau}$ does not depend on the second layer, the horizontal divergence is actually the standard Euclidean
divergence in $\mathbb{R}^{m}$. Therefore, using integration by parts, the following is
obtained:
\begin{align*}
V_{\delta, \tau}(\xi) &=\phantom{-}\int_{\mathbb{R}^{m}} V_{\delta, \tau} (\tilde{\xi})\,K_\tau(\xi-\tilde{\xi})\, d\tilde{\xi}, \textrm{ so,} \\
\textrm{div }V_{\delta, \tau} (\xi) &=\phantom{-}\int_{\mathbb{R}^{m}} \langle V_{\delta, \tau} (\tilde{\xi}),D_{\xi}(K_\tau(\xi-\tilde{\xi}))\rangle\, d\tilde{\xi}\\
&=-\int_{\mathbb{R}^{m}} \langle V_{\delta, \tau} (\tilde{\xi}),D_{\tilde{\xi}}(K_\tau(\xi-\tilde{\xi}))\rangle\, d\tilde{\xi}\\
&=-\int_{B(0,\delta)} \langle  -\frac{J(\epsilon_n)\tilde{\xi}}{\delta}  ,D_{\tilde{\xi}}(K_\tau(\xi-\tilde{\xi}))\rangle\, d\tilde{\xi}\\
&\phantom{=}-\int_{\mathbb{R}^{m} \backslash \overline{B(0,\delta)}} \langle
-\frac{J(\epsilon_n)\tilde{\xi}}{|\tilde{\xi}|}
,D_{\tilde{\xi}}(K_\tau(\xi-\tilde{\xi}))\rangle\, d\tilde{\xi}\\
&=\phantom{-}\int_{\partial B(0,\delta)} \langle -\frac{J(\epsilon_n)\tilde{\xi}}{\delta} , \nu_B
\rangle \,d\mathcal{H}^{m-1}\\
&\phantom{=}-\int_{\partial B(0,\delta)} \langle -\frac{J(\epsilon_n)\tilde{\xi}}{\delta} , \nu_B
\rangle \,d\mathcal{H}^{m-1} =0.
\end{align*}
Now applying the Gauss Green theorem with $V_{\delta, \tau}$ one obtains
\[
\int _{\mb E \cap \partial_{*, \gr} \Omega} \langle \nu _\Omega, V_{\delta, \tau} \rangle \, d\sddq =   \int
_{\mb E \cap \Omega} \langle \nu _E , V_{\delta, \tau} \rangle \, d\sddq.
\] 
One can easily check that the right hand side is less than or equal to $\sdq(\mb E~\cap~\Omega)$, for every $\tau >0$.  Since $\overline{\mb E \cap \mb \Omega}$ is a compact set, the
left hand side converges to 
\[
\int _{\mb E \cap \partial_{*, \gr} \Omega} \langle \nu _\Omega, V_{\delta} \rangle \, d\sddq
\]
as $\tau$ goes to $0$.  To compute this integral, it is split as
\begin{align*}
&\phantom{+}\int_{\mb E \cap \partial_{*, \gr} \Omega \cap \{ s > \delta\} } \langle -\frac{J(\epsilon_n)\xi}{s} , \nu_\Omega \rangle \, d\sdq\\
&+\int_{\mb E \cap \partial_{*, \gr} \Omega \cap \{
s \leq \delta\} } \langle -\frac{J(\epsilon_n)\xi}{\delta} , \nu_\Omega \rangle \, d\sdq\\
&= I_\delta + II_\delta.
\end{align*}
Using \ref{ht}, \ref{zero}, \ref{nmod}, and \ref{nvector}, one can compute
\begin{align*}
\langle -\frac{J(\epsilon_n)\xi}{s} , \nu_\Omega \rangle |\nabla_H \phi|&= \langle
\frac{J(\epsilon_n)\xi}{s}, \frac{\partial g}{\partial s} \frac{\xi}{s} \,\rangle + \langle
\frac{J(\epsilon_n)\xi}{s}, \frac{1}{2} J(\eta)\xi \,\rangle \\
&= 0 + \frac{1}{2s}\langle \epsilon_n, \eta \rangle |\xi|^2 \\
&= \frac{s}{2},
\end{align*}
and
\[
|\nabla_H \phi|^2 = \left(\frac{\partial g}{\partial s}\right) ^2 + \frac{1}{4}|\eta|^2|\xi|^2.
\]
Therefore,
\[
I_\delta= \int_{\mb E \cap \partial_{*, \gr} \Omega \cap \{ s > \delta\} }
\frac{1}{\sqrt{\left(\frac{\partial g}{\partial s}/s\right)^2 + \frac{1}{4}|\eta|^2}}\, d\sdq .
\]
Observe that the integrand is in $L^\infty(\mb E \cap \partial_{*, \gr} \Omega)$ by the hypothesis. In particular, there is $L(\Omega) >0$ such that $1 \leq |\eta|^2 = \sum_{l=1}^{n-1} \left(\frac{\partial
g}{\partial y_l}\right) ^2 +1 \leq L$ since $g$ is $C^1$ in the variables of the second layer.  Let $\delta$ go to $0$. By Lebesgue dominated convergence $I_\delta \to I$,
where
\[
I= \int_{\mb E \cap \partial_{*, \gr} \Omega}
\frac{1}{\sqrt{\left(\frac{\partial g}{\partial s}/s\right)^2 + \frac{1}{4}|\eta|^2}} \, d\sdq .
\]
Then by the hypothesis and above comment,
\[
I \geq \frac{1}{\sqrt{M^2 +\frac{1}{4}L}} \sdq (\mb E \cap \mb \Omega).
\]
Similarly, one may compute
\begin{align*}
II_\delta &= \frac{1}{\delta} \, \int_{\mb E \cap \mb \Omega \cap \{ s \leq \delta\}}
\frac{s^2 }{\sqrt{\left( \frac{\partial g}{\partial s}\right)^2 + \frac{1}{4}|\eta|^2 s^2} }\, d\sdq \\
&\leq  \int_{\mb E \cap \mb \Omega \cap \{ s \leq \delta\}} \frac{1}{\sqrt{\left(\frac{\partial g}{\partial s}/s\right)^2 + \frac{1}{4}|\eta|^2}} \, d\sdq .
\end{align*}
By Lebesgue dominated convergence, $II_\delta \to 0$ as $\delta \to 0$.  Combining all of
this information gives
\[
\sdq (\mb E \cap \mb \Omega ) \leq \sqrt{M^2 + \frac{1}{4}L} \ \,\sdq (\mb E \cap \Omega ).
\]
\end{proof}

The above theorem implies the $\gr$-admissibility of the gauge ball in the Heisenberg group. 

An example is presented showing that the condition $\Omega \in C^{1,\alpha}$, $\alpha <1$,
is not sufficient for $\gr$-admissibility. Therefore, the assumption $\Omega \in
C^{1,\alpha}$, $\alpha <1$ is not sufficient for $\gr$-admissibility for arbitrary Carnot
Groups of step 2. 

\begin{example}
Assume $\Omega \subset \mathbb{H}^1$, and $\partial \Omega$ is given by
$t=1-s^{2-\epsilon}$ in a neighborhood of
$P=(0,0,1)$. One can compute that
\[
\left| \frac{g'(s)}{s} \right| = \left| \frac {2-\epsilon}{s^\epsilon} \right|,
\] 
and observe that $g$ is not $C^{1,1}$ in a neighborhood of $P$.  
For every $B(P,1/N)$ a set $E_N \subset B(P,r) \cap \overline{\Omega}$ exists such that condition (ii)
of $\gr$-admissibility fails. 
In particular, 
\[
\frac{\sdq(\mb E_N \cap \mb \Omega)}{\sdq(\mb E_N \cap \Omega)}
\to \infty, \text{ as } N \to \infty.
\]
Consider $A_N := \{ t>1-1/N \}$, and let $E_N := A_N \cap \overline{\Omega}$. First, observe that $1-s^{2-\epsilon} =
1-1/N$ implies that $s=(1/N)^{\frac{1}{2-\epsilon}}$. A defining function $f$, for $A_N$ is
$f(s,t) = 1/N -t$. Therefore, by the structure theorem,
\[
 \sdq(\mb E_N \cap
\Omega)=\sdq(\mb A_N \cap
\Omega)=c_n \,|\partial A_N|_{\gr} (\Omega).
\]
By proposition \ref{form}, one computes
\begin{align*}
|\partial A_N|_{\gr} (\Omega) &=\int_{\mb A_N \cap \Omega} |\nabla_{\gr} f| \,
d\mathcal{H}^2\\
&=\int_{\sqrt{|x|^2+|y|^2}< (1/N)^{\frac{1}{2-\epsilon}}} |\nabla_{\gr}f| \, dx\,dy \notag \\ 
&= 2\pi \, \int_{0}^{ \left( \frac{1}{N} \right) ^ {\frac{1}{2-\epsilon}}}
\frac{r^2}{2} \, dr \notag \\
&=\frac{\pi}{3} \left( \frac{1}{N} \right) ^{\frac{3}{2-\epsilon}}. \notag
\end{align*}    
Next, observe that 
\[\sdq(\mb E_N \cap \mb \Omega)= \sdq(\mb \Omega \cap
A_N)= c_n\,|\partial \Omega|_\gr (A_N).
\]
The def{}ining
function for $\Omega$ is $g(s,t) = t+s^{2-\epsilon} -1$. Thus
\begin{align*}
|\partial \Omega|_{\gr} (A_N)&=\int_{\mb \Omega \cap A_N} |\nabla_{\gr} g| \,
d\mathcal{H}^2\\
&=\int_{\sqrt{|x|^2+|y|^2}< (1/N)^{\frac{1}{2-\epsilon}}}  |\nabla
_{\gr}g| \, dx\,dy \notag \\ 
&= 2\pi \, \int_{0}^{\left( \frac{1}{N} \right) ^ {\frac{1}{2-\epsilon}}}
 r \sqrt{(2-\epsilon)^2 r^{2(1-\epsilon)} + \frac{r^2}{4}}\, dr. \notag 
\end{align*}
Let 
\[
F(x) = \frac{1}{{x}^{\frac{3}{2-\epsilon}}} \, \int_{0}^{x^ {\frac{1}{2-\epsilon}}}
 r \sqrt{(2-\epsilon)^2 r^{2(1-\epsilon)} + \frac{r^2}{4}}\, dr.
\]
Then in order to show $\Omega$ is not $\gr$-admissible it suffices to show
\[
\lim_{x \to 0} F(x) = \infty.
\]
Using L'Hospital's rule, one observes
\begin{align}
\lim_{x \to 0} F(x) &= \lim_{x \to 0} \frac{1}{3} \,
\frac{1}{x^{\frac{1}{2-\epsilon}}}
\sqrt{(2-\epsilon)^2x^{\frac{2(1-\epsilon)}{2-\epsilon}} +
\frac{x^{\frac{2}{2-\epsilon}}}{4}} \notag \\
&= \lim_{x \to 0} \frac{1}{3} \, \sqrt{(2-\epsilon)^2
x^{\frac{-2\epsilon}{2-\epsilon}} + \frac{1}{4}} \notag \\
&= \infty. \notag
\end{align} 
\end{example}

\section{Approximate Continuity and Functions of $\gr$-Bounded Variation}
Now that the class of $\gr$-admissible domains has been shown to be non-empty, this section provides needed tools to prove an existence and trace theorem for $u \in \BVX(\Omega)$,
when $\Omega$ is $\gr$-admissible.

The upper approximate limit of $u$ at $x$, $\mu (x)$, and the approximate
lower limit of $u$ at $x$, $\lambda (x)$, are defined below:
\[
\mu (x) := \textrm{ap} \limsup _{y \rightarrow x} u(y) = \inf \{t:
D(A_t,x)=0 \},
\]
\[
\lambda (x) := \textrm{ap} \liminf _{y \rightarrow x} u(y) = \sup \{ t:
D(B_t, x) =0 \},
\]
where $A_t = \{x \in \gro : u(x) > t \}$, $B_t = \{ x \in \gro : u(x) < t \}$, and $D( \Omega,
x)$ is the density of $\Omega$ at $x$ with respect to the balls $B(x,r)$. The function
 $u$ is said to be approximately continuous at $x$ if $\lambda (x) = \mu (x)= u(x)$. 
The following proposition is proved as in remark 5.9.2 of \cite{Z}.

\begin{prop}\label{ac}
The function $u$ is approximately continuous at $x$ if and only if there exists
a Lebesgue measurable set $E$ containing $x$ such that $D(E,x) =1$ and $u\lfloor E$ is
continuous.
\end{prop}

A very significant result in the classical theory is that a function of bounded
variation is approximately continuous except for a set of $\mathcal{H}
^{n-1}$-measure zero, where it has a ``measure-theoretic jump,'' see \cite{Z}.  A
similar statement is true for functions of $\gr$-bounded variation with
Euclidean Hausdorff measure replaced with $\sdq$.  Throughout the rest of this
discussion, $F$ will denote the following:
\[
F= \{ x\in \gro: \lambda (x) < \mu (x) \}.
\]
\begin{thrm} \label{cool}
Let $u \in \BVX (\gro)$, then
\begin{itemize}
\item[(i)] $F$ is countably $\gr$-rectif{}iable,
\item[(ii)] $- \infty < \lambda (x) \leq \mu (x) < \infty$ for $\sdq$ a.e. $x \in \gro$,
\item[(iii)] for $\sdq$ a.e. $y \in F$, there is a  vector
$\nu(y)$ such that {} $\nu_{A_s}(y) = \nu(y)$ whenever $\lambda (y) < s < \mu
(y)$.
\item[(iv)] For all $y$ as in (iii), with $-\infty < \lambda (y) < \mu (y) <
\infty$ there are Lebesque measurable sets $F^+$ and $F^-$ such that
\begin{align*}
&\lim_{r \rightarrow 0} \frac{|F^- \cap \tau _y (S_\gr ^- ( \nu  (y))) \cap
B(y,r)|}{|\tau_y(S_\gr ^- ( \nu (y))) \cap B(y,r)|}\\
= &\lim_{r \rightarrow 0} \frac{|F^+ \cap \tau _y (S_\gr ^+( \nu  (y))) \cap
B(y,r)|}{|\tau _y (S_\gr ^+ ( \nu (y))) \cap B(y,r)|} =1
\end{align*}
and
\[
\lim_{\substack{x \rightarrow y \\ x \in F^- \cap \tau _y (S_\gr ^- ( \nu (y)))}}
u(x) = \mu (y), \qquad \lim_{\substack{x \rightarrow y \\ x \in F^+ \cap \tau _y
(S_\gr ^+ (\nu(y)))}}
u(x) = \lambda (y).
\]
\end{itemize}
\end{thrm}
\medskip
There are a few lemmas that will be used for the proof of the preceeding
theorem.
\begin{lemma} \label{per}
There exists $C=C(m+n)$ such that if $J$ is a Borel set such that $J \subset \rb
E$, then
\[
\sdq (J) \leq C \ \pme (J).
\]
\end{lemma}

\begin{proof}[Proof:]
Using lemma \ref{fss}, the proof of lemma 3.2.1 in \cite{Z} can be followed to reach the
conclusion.
\end{proof}
\begin{lemma} \label{vit}
Let $m+n \geq 1$ and $0 < \tau < 1/2$.  Let $E$ be a Lebesgue-measurable set such
that
\[
\lim_{ r \rightarrow 0} \frac{ |B(x,r) \cap E|}{|B(x,r)|} > \tau,
\quad \textrm{when }x \in E.
\]
Then there exists $C=C(\tau , m+n)$ and sequence of gauge balls $B_\rho (x_i, r_i)$ with $x_i \in
E$ such that
\[
E \subset \bigcup _{i=1}^\infty B_\rho (x_i, r_i),
\]
\[
\sum _{i=1} ^\infty (r_i) ^{Q-1} \leq C \ \pme (\gro).
\]
\end{lemma}

\begin{proof}

First observe that if the hypothesis of the lemma holds with $B(x,r)$ then it also holds with
the gauge ball $\br$, except with $\tau$ replaced with $\tau ' $, where $\tau '$  also
satisfies $0 < \tau ' < 1/2$.  To see this, first observe
that since there is $C>1$ such that
\[
C^{-1} \rho (x,y) \leq d(x,y) \leq C \ \rho (x,y),
\]  
 $\br \subset B(x,Cr)$.  There are constants $C_d >1$ and $C_\rho >1$ depending on the corresponding doubling constants of the respective
balls such that 
\[
|B_\rho (x,r)| \leq |B(x,Cr)| \leq C_d \ |B(x,r)| , \textrm{ and}
\]
\[
|B(x,r) \cap E | \leq |B_\rho(x,Cr)\cap E | \leq C_\rho |\br \cap E|.
\]
Therefore,
\[
\lim_{r \rightarrow 0}\frac{|B_\rho (x,r) \cap E|}{|B_\rho (x,r)|} \geq \lim_{r
\rightarrow 0} \frac{|B(x,r) \cap E|}{C_\rho \ C_d \ |B(x,r)|}
>\frac{\tau}{C_\rho \ C_d}.
\]
Then since $C_\rho, C_d$ are $> 1$,  $\tau ' := \frac{\tau}{C_\rho\ C_d}$
satisfies $0 < \tau ' < 1/2.$
Using theorem \ref{iso}, the proof can then
be completed as  lemma 5.9.3 in \cite{Z}. 
\end{proof}
\medskip
\begin{proof}[\textbf{Proof of theorem \ref{cool}:}]
The proof of theorem 5.9.6 in \cite{Z} is followed, but one must be extremely careful in applying
the various results stated thus far in this paper.  The proof in \cite{Z} makes use of the
set where the measure-theoretic normal exists.  This set can be replaced with the reduced
boundary without obstacle.  A proof of a
stronger statement which implies (i) can be found in
\cite{AM}, but a shorter proof is presented here. The complete proof of theorem \ref{cool} is given for the reader's benefit.

By the coarea formula for $\BVX$ functions, theorem \ref{coarea},  there 
exists a countable dense set $\mathcal{Q} \subset \mathbb{R}$ such that $|\partial A_t|_\gr
(\gro)
<\infty$ and $\rb A_t$ is countably $\gr$-rectifiable for $t \in \mathcal{Q}$.  By lemma
\ref{mea},
\[
\sdq [ \{ \cup ( \mb A_t \backslash  \rb A_t )  : t \in \mathcal{Q} \} ] =0.
\]
Using the definitions, one can observe
\begin{equation} \label{inc}
\{ x: \lambda (x) < t < \mu (x) \} \subset \mb A_t \textrm{  for } t \in
\mathbb{R}.
\end{equation}
Thus, $F\subset\{\cup \mb A_t :t\in\mathcal{Q}\}$, and $\sdq[F
\,\backslash\{ \cup
\rb A_t:t \in\mathcal{Q}\}]=0$.  Therefore, $F$ is countably $\gr$-rectifiable.  This
concludes the proof of (i).

For (ii), let $I=\{ x:\lambda (x)=\infty \} \cup\{x:\mu (x)=\infty \}$.
It be will shown that $\sdq (I) =0$.  Since this is a local question, it may be assumed that
 $u$ has compact support.  First it will be shown that $\sdq (\{x:\lambda (x)
= \infty \} )=0$ and $\sdq (\{ P: \mu (P) = -\infty \}) =0$ so that the set
$K=\{ x: \mu (x) - \lambda (x) = \infty\}$ is well-defined. Then it will be shown that
$\sdq (K) =0$, ending the proof of (ii).  By (i),  $u = \mu =
\lambda$ Lebesgue almost everywhere.  Therefore, letting $L_t = \{ x:\lambda (x)
> t\}$, 
\[
\frac{| \{ x: \lambda (x) > t \} \cap B(x, r) | }{|B(x,r)|} = \frac{ | \{ x: u(x)
>t \} \cap B(x,r) |}{ | B(x,r)|},
\]
and by definition, $\lambda (x) > t$ implies
\[
\lim_{ r \rightarrow 0} \frac{| \{x: u(x) < t \} \cap B(x, r)|}{|B(x,r)|} =0.
\]
Therefore, $D(L_t, x) =1$ for $x \in L_t$.  By lemma $\ref{vit}$, there exists a
countable family of gauge-balls, $\{ B_\rho (x_i, r_i) \}$ such that 
\[
L_t \subset \bigcup _{i=1}^\infty B_\rho (x_i, r_i) \subset \bigcup_{i=1}^\infty B(x_i,cr_i),
\]
\[
\sum _{i=1}^\infty (r_i) ^ {Q-1} \leq C |\partial L_t |_\gr (\gro).
\]
Since $u$ has compact support, one may assume $\diam \,B(x_i, cr_i) < a$.  From proposition \ref{coarea} and the fact $u= \lambda$ Lebesgue almost everywhere, 
\[
Var_\gr (u; \gro ) = Var_\gr (\lambda ; \gro) = \int_{-\infty}^\infty |\partial
L_t |_\gr (\gro) dt  < \infty, 
\]
implying
\[
 \lim _{t \rightarrow \infty} |\partial L_t| _\gr (\gro) = \lim _{t
\rightarrow \infty} \sdq (\mb L_t) =0.
\]
Thus,
\[
\mathcal{S}^{Q-1} _{d,a} [ x: \lambda (x) = \infty ] = \mathcal{S}^{Q-1} _{d,a} [ \bigcap _{ t=1 }^{\infty} L_t] \leq
C \sum _{i=1}^\infty (r_i)^{Q-1} \leq C \liminf_{t \rightarrow \infty }
|\partial L_t|_\gr (\gro ) =0.
\]
Let $a \rightarrow 0$.  Then $\sdq ( \{ x:\lambda (x)=\infty \} ) =0$.  Similarly,
$\sdq(\{x:\mu (x)=-\infty\}) =0.$  This concludes the proof of (ii).

From part (i),  $F$ is $\gr$-rectif{}iable.  Using the structure theorem
and the implicit function theorem, it can be concluded that $F$ is $\sigma$-finite
with respect to $\sdq \lfloor F$.  Thus,
\[
\begin{aligned}
\int _F ( \mu - \lambda )\, d \sdq  & = \int _0 ^\infty \sdq ( \{ x: \lambda (x) < t
< \mu (x) \} )\, dt \\
& \leq \int _0 ^\infty \sdq ( \mb A_t )\, dt  &\textrm{(by } \ref{inc} \textrm{)}\\
& \leq \int _0 ^\infty \sdq (\rb A_t ) \,dt &\textrm{(by lemma } \ref{mea} \textrm{)}\\
& \leq C \int _0 ^\infty |\partial A_t|_\gr (\gro )\, dt  &\textrm{(by lemma } \ref{per}\textrm{)}\\
& \leq C \ Var_H (u; \gro ) \\
& < \infty.
\end{aligned}
\]
Thus, $\sdq (K) =0.$

Part (iii) will be proved  for $x \in F \backslash \{ \cup ( \mb A_t
\backslash \rb A_t) : t \in \mathcal{Q} \}$. Without loss of generality, 
assume $x=0$.  One can check using the definitions that for $x \in \rb A_t, D(A_t, x) = 1/2$.
Now for $s < t, A_t \subset A_s$ and
\[
|A_s \backslash A_t \cap B(x,r)| = |A_s \cap B(x,r)| - |A_t \cap B(x,r)|
\]
implies that $D(A_s \backslash A_t, x) = 1/2 - 1/2 =0$.  This gives $\nu _{A_t} (0) = \nu _{A_s} (0)$ since
\[
\lim_{r \rightarrow 0} \frac{|B(x,r) \cap A_s \cap \mathcal{S}_\gr ^+ (\nu _{A_t} (0))|}{|B(x,r)|} 
\]
\[
= \lim_{r \rightarrow 0} \frac {|B(x,r) \cap A_t \cap
\mathcal{S}_\gr ^+ (\nu _{A_t} (0))|}{|B(x,r)|} + \lim _{r \rightarrow 0}
\frac{|B(x,r) \cap (A_s \backslash A_t) \cap \mathcal{S}_\gr ^+ (\nu _{A_t}
(0))|}{|B(x,r)|} =0, 
\]
and 
\begin{align*}
\lim _{r \rightarrow 0} \frac{|(B(x,r) \backslash A_s) \cap \mathcal{S}_\gr ^-
(\nu _{A_t} (0))|}{|B(x,r)|}&= \lim_{r \rightarrow 0}\frac{|B(x,r) \cap \{ x: u < s \} \cap
\mathcal{S}_\gr ^- (\nu _{A_t}
(0))|}{|B(x,r)|}\\
&\leq \lim _{r \rightarrow 0}\frac{|B(x,r) \cap \{ x: u < t \} \cap \mathcal{S}_\gr ^- (\nu _{A_t}
(0))|}{|B(x,r)|}\\
&= \lim _{r \rightarrow 0} \frac{|(B(x,r) \backslash A_t) \cap \mathcal{S}_\gr ^- (\nu _{A_t}
(0))|}{|B(x,r)|} =0.
\end{align*}
This concludes the proof of (iii).

In proving (iv), consider $y \in F \backslash I$.  Let $\epsilon > 0$ be such
that $\lambda (Q)<\mu (Q)-\epsilon < \mu (Q)$.  By definition, $D( A_{\mu
(y) + \epsilon} , y)=0$, and one can observe that
\[
\lim _{r \rightarrow 0} \frac{|A_{\mu (y) -
\epsilon} \cap \mathcal{S}_\gr ^- (\nu _{A_{\mu(y) -\epsilon}}
(0))\cap B(y,r)|}{|B(y,r) \cap \mathcal{S}_\gr ^- (\nu _{A_{\mu(y) -\epsilon}}
(0))|} = 1.
\]
Further, $A_{ \mu (y) - \epsilon } \subset u^{-1} \{ \lbrack \mu (y) - \epsilon, \mu
(y) + \epsilon \rbrack \}$. Therefore,
\medskip
\[
\lim _{r \rightarrow 0} \frac{|u^{-1} \{ \lbrack \mu (y) - \epsilon, \mu
(y) + \epsilon \rbrack \} \cap  \mathcal{S}_\gr ^- (\nu _{A_{\mu(y) -\epsilon}}
(0)) \cap B(x,r)|}{|B(x,r) \cap  \mathcal{S}_\gr ^- (\nu _{A_{\mu(y) -\epsilon}}
(0))|} =1.
\]
Similarly to proposition \ref{ac}, this gives  the desired set $F^-$.  Similarly, 
$F^+$ is obtained.  
\end{proof}

\section{Defining the trace of an $\gr$-BV function on an $\gr$-admissible domain}

This section proves an extension and trace theorem for $u \in \BVX(\Omega)$, when $\Omega$ is $\gr$-admissible. 

\begin{dfn}
Let $u$ be a real-valued Lebesgue measurable function def{}ined on $\Omega$, an
open set in $\gro$. 
\[
u_0 (x) := \begin{cases}
u(x)& \text{if $x \in \Omega$},\\
0& \text{if $x \in \Omega ^C$}.
\end{cases}
\]
\end{dfn}
\medskip
\noindent \textbf{Observation:  } Let $F= \{ x : \lambda_{u_0} < \mu_{u_0} \}$, $A_t= \{ x:
u_0 > t \}$, and $x_0 \in E \cap \rb \Omega$ be such that (iii) of Theorem \ref{cool} applies.  Then
there is a unit vector $\nu$ such that
\[
\nu_{A_t}(x_0) = \nu \qquad \textrm{whenever } \lambda_{u_0} < t < \mu_{u_0}.
\]
As in Remark 5.10.6 in \cite{Z}, one can conclude that
\begin{equation}\label{nvec}
\nu = \pm \nu_\Omega (x_0),
\end{equation}
and that 
\begin{equation}\label{posn}
\textrm{if } \nu = \nu_{\Omega}(x_0), \textrm{ then } \lambda_{u_0} (x_0) =0,
\end{equation}
\begin{equation}\label{negn}
\textrm{if } \nu = -\nu_{\Omega}(x_0), \textrm{ then } \mu_{u_0} (x_0) =0.
\end{equation}

\begin{prop}
If $\Omega$ is an $\gr$-admissible domain and $u \in \BVX(\Omega )$, then $u_0 \in
\BVX(\gro)$ and there exist $C=C(\Omega)$ such that
\[
\parallel u_0 \parallel _{\BVX (\gro)} \leq C \ \parallel u \parallel _{\BVX
(\Omega)}.
\]
\end{prop}
\medskip
\begin{proof}[Proof:]
It is a fact that $u \in \BVX (\Omega)$ if and only if the positive and negative parts of
$u$, $u^+$ and $u^-$, belong to $\BVX (\Omega)$.  This fact along with the coarea formula for
functions in $\BVX$ and the structure theorem provide the proof for the theorem analogously
to lemma 5.10.4 in \cite{Z}.
\end{proof}
\medskip
Now  the trace of a function on $\partial \Omega$ for $\Omega$ an
$\gr$-admissible domain is given:
\begin{dfn}
If $\Omega$ is an $\gr$-admissible domain and $u \in \BVX (\Omega)$, the trace of $u$ on
$\partial \Omega$, denoted
$u^*$,  is
\[
u^* (x) := \mu _{u_0} (x) + \lambda _{u_0} (x),
\]
where, $\mu _{u_0} , \lambda _{u_0}$ are the upper and lower approximate limits of
$u_0$.
\end{dfn}
\medskip
\begin{thrm}
If $\Omega$ is an $\gr$-admissible domain, there is a constant $M = M(\Omega)$
such that
\[
\int _{\gro} |u^*| \ \ d |\partial \Omega|_\gr = C(n+m) \ \int _{\mb \Omega} |u^*| \
\ d\sdq \leq  M \ \parallel u \parallel
_{\BVX (\Omega)}
\]
whenever $u \in \BVX (\Omega).$
\end{thrm}

\begin{proof}[Proof: ]
The proof follows as theorem 5.10.7 in \cite{Z} where the coarea formula for $\BVX$ is
utilized, as is theorem \ref{cool} (i).
\end{proof}
\medskip

\section{The trace and integral averages}
Throughout this section $F$ is defined as in previous sections and the function $U$ is
defined as $U(x) = (\lambda_u (x) + \mu_u(x) )/2$.  The following theorem will be proved:  
\begin{thrm}\label{mt}
Assume $u \in \BVX (\gro)$.  Then 
\begin{itemize}
\item[(i)] ${\lim_{r \rightarrow 0} \int \!\!\!\!\! - }_{\br}|u-U(x)| ^{\frac{Q}{Q-1}} \,dh =0$  for 
$\sdq$  a.e.  $x \in
\gro\backslash F$, and
\item[(ii)]for $\sdq$ a.e. $x \in F$, there exists a  vector $\nu =
\nu(x)$ such that
\[
\lim_{r \to 0}  {\int \!\!\!\!\!\! - }_{\!\!\br \cap S_\gr^- (\nu)} |u - \mu (x)|
^\frac{Q}{Q-1} \, dh = 0,
\]
and
\[
\lim_{r \to 0}  {\int \!\!\!\!\!\! - }_{\!\!\br \cap S_\gr^+ (\nu)} |u - \lambda (x)|
^\frac{Q}{Q-1} \, dh =0.
\]
\end{itemize}
\end{thrm}

The following corollary follows from \ref{nvec}, \ref{posn}, and \ref{negn}:
\begin{cor}
Let $u \in \BVX(\Omega)$, where $\Omega$ is $\gr$-admissible.  Then for $\sdq$ a.e. $x_0 \in
\partial \Omega$,
\[
\lim_{r \to 0} \int_{B(x_0, r) \cap \Omega} |u(x) - u^*(x_0)|^{\frac{Q}{Q-1}} \, dh =0.
\]
\end{cor}

The proof of theorem \ref{mt} relies on the following proposition:
\begin{prop}\label{in}
For every $0 < \alpha \leq 1$, there exists $C (\alpha)$ such that
\[
\|u\|_{L^\frac{Q}{Q-1}(\br)} \leq C \ Var_\gr (u ; \br)
\]
for all $\br \subset\gro$ and all $u \in BV_{\gr, loc} (\gro)$ such that
\[
\frac{|\br \cap \{ x: u(x) =0 \}|}{|\br|} \geq \alpha.
\]
\end{prop}
  
\begin{proof}
Assume that
\[
\frac{|\br \cap \{ u=0 \} |}{|\br|} \geq \alpha > 0.
\]
Then,
\begin{align*}
\|u\| _{L^\frac{Q}{Q-1}(\br} & \leq  \|u -
u_{x,r}\|_{L^\frac{Q}{Q-1}(\br)} + \|u_{x,r}\| _{L^\frac{Q}{Q-1}(\br)}\\
& \leq  C \,Var_\gr(u; \br) + |u_{x,r}| \, |\br| ^{1-1/Q}.
\end{align*}
The last inequality follows from theorem \ref{poin}.  Next observe,
\begin{align*}
|u_{x,r}| \, |\br| ^{1-1/Q} & = \frac{|\br|}{|\br)| ^{1/Q}} \,\iavgr
|u| \, dh\\
& =  \frac{1}{|\br| ^{1/Q}} \int_{\br \cap \{u \neq 0\}} |u| \, dh\\
& \leq  \left( \int_{\br} |u| ^{\frac{Q}{Q-1}} \, dh \right) ^ {1-1/Q} \, \left(\frac{|\br
\cap \{ u \neq 0 \} |}{|\br|} \right) ^{1/Q} \\
& \leq \|u\| _{L^\frac{Q}{Q-1}(\br)} (1-\alpha) ^{1/Q}.
\end{align*}
The conclusion follows.
\end{proof}

\begin{proof}[\textbf{Proof of theorem \ref{mt}:}]
First, observe that replacing $B(x,r)$ by $\br$ in the definition of $\mu_u(x)$ and
$\lambda_u(x)$ result in the same function, since $\rho$ and $d$ are equivalent metrics.  The proof
now follows as theorem 3, section 5.9 in \cite{EG}. In this proof, theorem \ref{cool} (iii) and
proposition \ref{in} are vital.
\end{proof} 

\bibliographystyle{amsplain}
\bibliography{bibfile}

\def\cprime{$'$}
\providecommand{\bysame}{\leavevmode\hbox to3em{\hrulefill}\thinspace}
\providecommand{\MR}{\relax\ifhmode\unskip\space\fi MR }
\providecommand{\MRhref}[2]{%
  \href{http://www.ams.org/mathscinet-getitem?mr=#1}{#2}
}
\providecommand{\href}[2]{#2}
\begin{thebibliography}{10}

\bibitem{AMB}
Luigi Ambrosio, \emph{Fine properties of sets of finite perimeter in doubling
  metric measure spaces}, Set-Valued Anal. \textbf{10} (2002), no.~2-3,
  111--128, Calculus of variations, nonsmooth analysis and related topics.
  \MR{MR1926376 (2003i:28002)}

\bibitem{A}
Luigi Ambrosio, Nicola Fusco, and Diego Pallara, \emph{Functions of bounded
  variation and free discontinuity problems}, Oxford Mathematical Monographs,
  The Clarendon Press Oxford University Press, New York, 2000. \MR{MR1857292
  (2003a:49002)}

\bibitem{AM}
Luigi Ambrosio and Valentino Magnani, \emph{Weak differentiability of {BV}
  functions on stratified groups}, Math. Z. \textbf{245} (2003), no.~1,
  123--153. \MR{MR2023957 (2005a:43010)}

\bibitem{BA}
Zolt{\'a}n~M. Balogh, \emph{Size of characteristic sets and functions with
  prescribed gradient}, J. Reine Angew. Math. \textbf{564} (2003), 63--83.
  \MR{MR2021034 (2005d:43007)}

\bibitem{CDG1}
Luca Capogna, Donatella Danielli, and Nicola Garofalo, \emph{The geometric
  {S}obolev embedding for vector fields and the isoperimetric inequality},
  Comm. Anal. Geom. \textbf{2} (1994), no.~2, 203--215. \MR{MR1312686
  (96d:46032)}

\bibitem{CG}
Luca Capogna and Nicola Garofalo, \emph{Boundary behavior of nonnegative
  solutions of subelliptic equations in {NTA} domains for
  {C}arnot-{C}arath\'eodory metrics}, J. Fourier Anal. Appl. \textbf{4} (1998),
  no.~4-5, 403--432. \MR{MR1658616 (2000k:35056)}

\bibitem{DGN1}
Donatella Danielli, Nicola Garofalo, and Duy-Minh Nheiu, \emph{Minimal surfaces
  in carnot groups}, preprint (2004).

\bibitem{EG}
Lawrence~C. Evans and Ronald~F. Gariepy, \emph{Measure theory and fine
  properties of functions}, Studies in Advanced Mathematics, CRC Press, Boca
  Raton, FL, 1992. \MR{MR1158660 (93f:28001)}

\bibitem{FE}
Herbert Federer, \emph{Geometric measure theory}, Die Grundlehren der
  mathematischen Wissenschaften, Band 153, Springer-Verlag New York Inc., New
  York, 1969. \MR{MR0257325 (41 \#1976)}

\bibitem{FSSC3}
Bruno Franchi, Raul Serapioni, and Francesco Serra~Cassano,
  \emph{Meyers-{S}errin type theorems and relaxation of variational integrals
  depending on vector fields}, Houston J. Math. \textbf{22} (1996), no.~4,
  859--890. \MR{MR1437714 (98c:49037)}

\bibitem{FSSC2}
\bysame, \emph{Rectifiability and perimeter in the {H}eisenberg group}, Math.
  Ann. \textbf{321} (2001), no.~3, 479--531. \MR{MR1871966 (2003g:49062)}

\bibitem{FSSC1}
\bysame, \emph{On the structure of finite perimeter sets in step 2 {C}arnot
  groups}, J. Geom. Anal. \textbf{13} (2003), no.~3, 421--466. \MR{MR1984849
  (2004i:49085)}

\bibitem{GN}
Nicola Garofalo and Duy-Minh Nhieu, \emph{Isoperimetric and {S}obolev
  inequalities for {C}arnot-{C}arath\'eodory spaces and the existence of
  minimal surfaces}, Comm. Pure Appl. Math. \textbf{49} (1996), no.~10,
  1081--1144. \MR{MR1404326 (97i:58032)}

\bibitem{MAG}
Valentino Magnani, \emph{Characteristic points, rectifiability and perimeter
  measure on stratified groups}, preprint (2003).

\bibitem{MZ}
Norman~G. Meyers and William~P. Ziemer, \emph{Integral inequalities of
  {P}oincar\'e and {W}irtinger type for {BV} functions}, Amer. J. Math.
  \textbf{99} (1977), no.~6, 1345--1360. \MR{MR0507433 (58 \#22443)}

\bibitem{Z}
William~P. Ziemer, \emph{Weakly differentiable functions}, Graduate Texts in
  Mathematics, vol. 120, Springer-Verlag, New York, 1989, Sobolev spaces and
  functions of bounded variation. \MR{MR1014685 (91e:46046)}

\end{thebibliography}
\end{document}